\documentclass[10pt]{article}
\usepackage{amsmath}
\usepackage{color}
\usepackage{amssymb}
\usepackage{epsfig}
\numberwithin{equation}{section}
\newcounter{indiceC}
\newcounter{indiceK}
\newcommand{\R}{\mathbb{R}}

\newcommand{\Z}{\mathbb{Z}}
\newcommand{\T}{\mathbb{T}}
\newcommand{\be}{\begin{equation}}
\newcommand{\ee}{\end{equation}}
\newcommand{\ba}{\begin{eqnarray}}
\newcommand{\ea}{\end{eqnarray}}

\newcommand{\cqfd}
{%
\mbox{}%
\nolinebreak%
\hfill%
\rule{2mm}{2mm}%
\medbreak%
\par%
}

\newtheorem{thm}{Theorem}

\newtheorem{rmk}{Remark}

\begin{document}


\title{On the Benjamin-Bona-Mahony equation with a localized damping}

\author{Lionel Rosier \thanks{Centre Automatique et Syst\`emes,
MINES ParisTech, PSL Research University, 60 Boulevard Saint-Michel, 
75272 Paris Cedex 06, France
Email: {\tt Lionel.Rosier@mines-paristech.fr}.}}

\date{\empty}
\maketitle

\begin{center}
{\bf Abstract}
\end{center}
We introduce several mechanisms to dissipate the energy in the Benjamin-Bona-Mahony (BBM) equation. We consider
either a distributed (localized) feedback law, or a boundary feedback law. 
In each case, 
we prove the global wellposedness of the system and the convergence towards a solution of the BBM equation which is null on 
a band. If the Unique Continuation Property holds for the BBM equation, this implies that the origin is 
asymptotically stable for the damped BBM equation. 


\vskip0.3cm\noindent {\it AMS Subject Classification:} 35Q53,  93D15. 
\vskip0.3cm\noindent {\it Keywords:} Benjamin-Bona-Mahony equation, unique
continuation property, internal stabilization, boundary stabilization.

\section{Introduction}\label{sec1}

The Benjamin-Bona-Mahony equation 
\be
\label{A1}
v_t-v_{xxt}+v_x+vv_x=0,
\ee
was proposed in \cite{BBM} as an alternative to the Korteweg-de Vries (KdV) equation as a model for the propagation of
one-dimensional, unidirectional, small amplitude long waves in nonlinear dispersive media. In the context of shallow water waves, 
$v=v(x,t)$ stands for the displacement of the water surface (from rest position) at location $x$ and time $t$. 
In the paper, we shall assume that either $x\in \R$, or $x\in (0,L)$ or $x\in \T=\R /(2\pi \Z )$ (the one-dimensional torus).

The dispersive term $-v_{xxt}$ produces a strong smoothing effect for the time regularity, thanks to which the wellposedness
theory of \eqref{A1} is easier than for KdV (see \cite{BT,roumegoux}). Solutions of \eqref{A1} turn out to be analytic in time.
On the other hand, the control theory is at his early stage for BBM (for the control properties of KdV, we refer the reader to 
the recent survey \cite{RZ2009}). S. Micu investigated in \cite{micu}
the boundary controllability of the {\em linearized} BBM equation, and noticed that the exact controllability fails due to the existence of a {\em limit point} in the spectrum of the adjoint equation. 
The author and B.-Y. Zhang  introduced in \cite{RZ2013} a {\em moving control}  and derived with such a control both the exact
controllability and the exponential stability of the full BBM equation.  For a distributed control with a {\em fixed support}, 
the exact controllability of the linearized BBM equation fails, so that the study of the controllability of the full BBM equation seems hard.
However, it is reasonable to expect that some stability properties could be derived by incorporating some dissipation  
in a fixed subdomain or at the boundary. The aim of this paper is to propose several dissipation mechanisms leading to systems
for which one has both the global existence of solutions and a  nonincreasing $H^1$-norm. 
The conclusion that all the trajectories emanating from data close to the origin are indeed attracted by the origin
is valid provided that the following conjecture is true: \\[3mm]
{\bf Unique Continuation Property (UCP) Conjecture:} There exists some number $\delta >0$ such that for any 
$v_0\in H^1(\T)$ with $\Vert v_0\Vert _{H^1(\T )} <\delta$, if the solution 
$v=v(x,t)$ of 
\be 
\label{A2}
\left\{ 
\begin{array}{ll}
v_t-v_{xxt}+v_x+vv_x =0 , \quad &x\in \T ,\\
v(x,0)=v_0(x), \quad &x\in\T
\end{array}
\right.
\ee
satisfies 
\be
\label{A3}
v(x,t)=0 \, \quad \forall (x,t)\in \omega \times (0,T) 
\ee
for some nonempty open set $\omega \subset \T$ and some time $T>0$, then $v_0=0$ (and hence $v\equiv 0$).  

To the best knowledge of the author, the UCP for BBM as stated in the above conjecture is still open. The main difficulty
comes from the fact that the lines $x=0$ are {\em characteristic} for BBM, so that the ``information'' does not propagate well 
in the $x$-direction. For some UCP for BBM  (with additional assumptions) see \cite{RZ2013,yamamoto}.   See also \cite{MORZ, MP} for 
control results for some Boussinesq systems of BBM-BBM type.
 
 The following result is a {\em conditional} UCP in which it is assumed that the initial data is small in 
 the $L^\infty$-norm and it has a {\em nonnegative} mean value.  Its proof was based on the analyticity in time of the trajectories and 
 on the use of some Lyapunov function. 

\begin{thm} \cite{RZ2013}
Let $u_0\in H^1(\T ) $ be such that 
\be
\label{AA10}
\int_\T v_0(x)\, dx \ge 0 , 
\ee
and 
\be
\label{AA11}
\Vert v_0\Vert _{L^\infty (\T ) } <3.
\ee
Assume that the solution $v$ to \eqref{A2} satisfies \eqref{A3}. Then $v_0=0$. 
\end{thm}

As it was noticed in \cite{RZ2013}, the UCP for BBM cannot hold for any state  in $L^2(\T)$, for any initial data $v_0$ with values in 
$\{ -2, 0 \}$ gives a trivial (stationary) solution of BBM. Thus, either a regularity assumption ($v_0\in H^1(\T)$), or a bound on the norm 
of the initial data has to be imposed for the UCP to hold. 

The paper is outlined as follows. In Section 2, we incorporate a simple localized damping in BBM equation and investigate 
the corresponding Cauchy problem. In Section 3, we consider another dissipation mechanism involving one derivative. The last
section is concerned with the introduction of nonhomogeneous boundary conditions leading again to the dissipation of the $H^1$-norm.

\section{Stabilization of the BBM equation}
\subsection{Internal stabilization with a simple feedback law}
We consider the BBM equation on $\T$ with a localized damping 
\begin{eqnarray}
u_t-u_{xxt}+u_x+uu_x +a(x)u =0 \quad && x\in \T,\ t\ge 0,\label{A20}\\
u(x,0)=u_0(x) \quad && x\in \T , \label{A21}
\end{eqnarray}
where $a$ is a smooth, nonnegative function on $\T$ with $a(x)>0$ on a given 
open set $\omega \subset \T$. 

We write \eqref{A20}-\eqref{A21} in its integral form
\begin{equation}
\label{A22}
u(t)=u_0 - \int_0^t (1-\partial _x^2)^{-1} 
[a(x)u + ( u + \frac{u^2}{2} )_x ] (\tau )d\tau
\end{equation}
where $(1-\partial _x^2)^{-1}f$  denotes for $f\in L^2(\T )$ the unique solution $v\in H^2(\T )$ of $(1-\partial _x^2)v=f$. 
Define the norm $\Vert \cdot \Vert _s$ in $H^s(\T )$ as 
\[
\Vert \sum_{n\in \Z} c_ne^{inx} \Vert _s^2 = \sum_{n\in \Z} \vert c_n\vert ^2 (1+\vert n\vert ^2)^s.
\]
We have the following result.
\begin{thm}
\label{stab1} 
Let $s\ge 0$ be given. For any $u_0\in H^s(\T )$, there exist
$T>0$ and a unique solution $u$ of \eqref{A22} in 
$C([0,T],H^s(\T ))$. Moreover, the correspondence $u_0\in H^s(\T )\mapsto 
u\in C([0,T],H^s(\T ))$ is Lipschitz continuous. If $s=1$, the solution exists for every $T>0$, and the energy
$\Vert u ( t ) \Vert _{1}  ^2$ is nonincreasing. Finally, 
if $v_0$ denotes any weak limit in $H^1(\T )$ of a sequence 
$u(t_n)$ with $t_n\nearrow +\infty$, then the solution $v$ of system
\eqref{A2} satisfies \eqref{A3}. In particular, if  the UCP conjecture holds,
then $v_0=0$, so that  $u(t)\to 0$ weakly in $H^1(\T)$ as $t\to +\infty$, hence 
strongly in  $H^s(\T )$ for all $s<1$.
\end{thm}
\noindent
{\em Proof.} We proceed as in \cite{BCS2}  using the estimate
\begin{equation}
\label{A15}
||fg||_s \le  C ||f||_{s+1} ||g||_{s+1} \qquad 
\forall s\ge -1,\ \forall f,g \in H^{s+1} (\T ).
\end{equation}
The estimate \eqref{A15} follows from a similar estimate proved in  \cite{BCS2} for functions defined on $\R$, namely
\begin{equation}
\label{A31}
||\tilde f \tilde g ||_{H^s(\R )} 
\le C ||\tilde f  ||_{H^{s+1}(\R )} ||\tilde g  ||_{H^{s+1}(\R )}\qquad
 \forall s\ge -1,\ \forall f,g \in H^{s+1}(\R ),
\end{equation}
by letting $\tilde f(x) =\varphi (x) f(x)$, $\tilde g(x)=\varphi (x)g(x)$, 
where $f$ and $g$ are viewed as $2\pi$-periodic functions, and $\varphi
\in C^\infty_0(\R )$ denotes a cut-off function such that $\varphi (x)=1$ on 
$[0,2\pi ]$. Indeed, we notice that 
\[
||f||_{H^s(\T )} \le ||\tilde f||_{H^s(\R )} \le C || f ||_{H^s(\T )}
\]
for some constant $C>0$. 
Note that for any $s \ge 0$ 
\[
||(1-\partial _x^2)^{-1} \partial _x (fg)||_s \le C ||fg||_{s-1} 
\le C ||f||_s||g||_s.
\]
Pick any $u_0\in H^s (\T )$. Let us introduce the operator
\[
 (\Gamma u)(t) := u_0  -  \int_0^t (1-\partial _x ^2)^{-1}
[a(x)u + (u + \frac{u^2}{2})_x ] (\tau ) d\tau .
\]
Then 
\begin{eqnarray*}
\sup_{0\le t\le T } ||(\Gamma u-\Gamma v)(t)||_s 
&\le & C_2 \int_0^T [||a(u-v)||_{s-2} + ||u-v||_{s-1}+||u-v||_s||u+v||_s] d\tau \\
&\le & C_3T (1+2R) ||u-v||_{C([0,T], H^s(\T ))} 
\end{eqnarray*}
if we assume that $u$ and $v$ are in the closed ball $B_R$ of radius $R$
centered at $0$ in $C([0,T],H^s(\T ))$. We pick 
$T:=[ 2C_3(1+2R)]^{-1}$ so that 
\[
||\Gamma u - \Gamma v||_{C([0,T],H^s(\T ))} \le \frac{1}{2} ||u-v||_{C([0,T],H^s(\T ))}\cdot
\]
On the other hand
\[
||\Gamma u||_{C([0,T],H^s(\T ))} 
\le ||u_0||_s + ||\Gamma u -\Gamma 0||_{C([0,T],H^s(\T ))} \le ||u_0||_s +\frac{R}{2}\le R
\]
for the choice $R=2||u_0||_s$. It follows that the map $\Gamma $ contracts in $B_R$, hence it 
admits a unique fixed point $u$ in $B_R$ which solves the integral equation 
\eqref{A22}. Furthermore, given any $\rho >0$ and any 
$u_0,v_0\in H^s(\T )$ with $||u_0||_s\le \rho$, $||v_0||_s\le \rho$, one easily sees
that for 
\begin{equation}
\label{A30}
T=[2C_3(1+4\rho )]^{-1}
\end{equation}
one has
\begin{equation}
\label{ABC}
||u-v||_{C([0,T],H^s(\T ))} \le 2||u_0-v_0||_s.
\end{equation}
Finally assume that $s=1$. Scaling in \eqref{A20} by $u$ yields
\begin{equation}
\label{A32}
\frac{1}{2}||u(T)||^2_1 -\frac{1}{2}||u_0||^2_1 
+\int_0^T\!\!\!\int_{\T }a(x)|u(x,t)|^2dxdt =0 
\end{equation}
(Note that $u_t= -(1-\partial _x^2)^{-1}[a(x)u  + (u+\frac{u^2}{2})_x]
\in C([0,T], H^2(\T ))$ so that each term in \eqref{A20} belongs to 
$C([0,T],L^2(\T ))$, and the integrations by parts are valid.) 
It follows that the map $t\mapsto ||u(t) || _1$ is nonincreasing, hence it admits
a nonnegative limit $l$ as $t\to \infty$, and that the solution 
of \eqref{A20}-\eqref{A21} is defined for all $t\ge 0$. Let $T$ be as in 
\eqref{A30} with $\rho = ||u_0||_1$. Note that 
$||u(t)||_1\le ||u_0||_1$ for all $t\ge 0$. Let $v_0$ be any weak limit 
of $\{ u(t) \} _{t\ge 0}$ in $H^1(\T )$, i.e. we have that
$u(t_n)\to v_0$ weakly in $H^1(\T )$ for some sequence $t_n\to + \infty$. Extracting a subsequence if needed, we
may assume that 
\begin{equation}
\label{Atn}
t_{n+1}-t_n\ge T\qquad \text{ for } n\ge 0.
\end{equation}
\null From 
\[
\frac{1}{2} ||u(t_{n+1})||_1^2 - \frac{1}{2} ||u(t_n)||_1^2
+\int_{t_n}^{t_{n+1}} \!\!\!\int_{\T} a(x)|u(x,t)|^2dxdt =0,
\] 
we obtain that 
\begin{equation}
\label{A35}
\lim_{n\to +\infty} \int_{t_n}^{t_{n+1}} \!\!\! \int_T  a(x)|u(x,t)|^2 dxdt=0.
\end{equation}
Since $u(t_n)\to v_0$ in $H^s(\T )$ for $s<1$, and 
$||u(t_n)||_s \le ||u(t_n)||_1 \le \rho$, we have from \eqref{ABC} that 
for all $s\in [0,1[$ 
\begin{equation}
\label{A36}
u(t_n+\cdot ) \to v\qquad \text{ in } C([0,T],H^s(\T ))\qquad
\text{ as } n\to +\infty , 
\end{equation}
where $v=v(x,t)$ denotes the solution of 
\begin{eqnarray*}
&& v_t - v_{xxt} + v_x + vv_x + a(x) v =0,\qquad x\in \T, \ t\ge 0, \\
&& v(x,0)=v_0(x),\qquad x\in \T .
\end{eqnarray*}
Notice that $v\in C([0,T];H^1(\T )) $, for $v_0\in H^1(\T )$.
\eqref{A35} combined to \eqref{Atn} and \eqref{A36} yields
\begin{equation}
\label{A37}
\int_0^T\!\!\!\int_{\T }a(x)|v|^2 dxdt=0.
\end{equation}
Therefore $v\in C([0,T];H^1(\T ))$ solves 
\begin{eqnarray*}
&&v_t-v_{xxt}+v_x+vv_x =0,\qquad x\in \T,\ t\in (0,T),\\
&&v(x,0)=v_0(x) \qquad x\in \T ,\\
&&v(x,t)=0 \qquad x\in \omega,\ t\in (0,T).
\end{eqnarray*}
If the UCP conjecture is true, we have that $v\equiv 0$ on $\R \times (0,T)$. It follows that 
$v_0\equiv 0$, and that as $t\to + \infty$ 
\begin{eqnarray*}
&&u(t)\to 0 \qquad \text{ weakly   in } H^1(\T ),\\
&&u(t)\to 0 \qquad \text{ strongly in } H^s(\T ) \text{ for any } s<1.
\end{eqnarray*}
\cqfd
\subsection{Internal stabilization with one derivative in the feedback law}
We now pay attention to the stabilization of the BBM equation by means of a 
``stronger'' feedback law involving one derivative. More precisely, we consider the system 
\begin{eqnarray}
&&u_t-u_{xxt}+u_x+uu_x -(a(x)u_x)_x =0, \qquad x\in \T, \ t\ge 0, \label{A50}\\
&&u(x,0) = u_0(x), \qquad x\in \T. \label{A51}
\end{eqnarray}
where $a=a(x)$ denotes again any nonnegative smooth function on $\T$ 
with $a(x)>0$ on a given open set $\omega \subset \T$.
Scaling in \eqref{A50} yields (at least formally)
\begin{equation}
\label{A52}
\frac{1}{2} ||u(T)||^2_1 -\frac{1}{2}||u_0||_1^2
+\int_0^T\!\!\!\int_{\T} a(x)|u_x(x,t)|^2 dxdt =0.
\end{equation}
The decay of the energy is quantified by an integral term involving
the square of a localized $H^1$-norm in \eqref{A52}. By contrast, the integral term in  \eqref{A32} 
involved  the square of a localized $L^2$-norm.   This suggests that the damping mechanism
involved in \eqref{A50} acts in a {\em much stronger way} than in \eqref{A20}. 
As a matter of fact, in the trivial situation when $a(x)\ge C>0$ for all $x\in \T$
in \eqref{A50}, it is a simple exercise to establish the exponential stability in $H^1(\T)$ for 
both the linearized equation and the nonlinear BBM equation for states with zero means.
In the general case when the function $a(x)$ is supported in a subdomain of $\T$, we obtain the following result. 
\begin{thm}
\label{stab2} 
Let $s\ge 0$ be given. For any $u_0\in H^s(\T )$, there exist
$T>0$ and a unique solution $u$ of \eqref{A50}-\eqref{A51} in 
$C([0,T],H^s(\T ))$. Moreover, the correspondence $u_0\in H^s(\T )\mapsto 
u\in C([0,T],H^s(\T ))$ is Lipschitz continuous. If $s=1$, one can pick any
$T>0$  and $\Vert u(t) -[u_0]\Vert _{H^1(\T )}  $ is nondecreasing. If the UCP conjecture is valid 
by replacing \eqref{A3} by $v=C$ on $\omega \times (0,T)$, then $u(t)\to [u_0]=(2\pi)^{-1} \int _\T u_0(x)dx$ weakly in $H^1(T)$  
hence  strongly in  $H^s(\T )$ for $s<1$ as $t\to \infty$.
\end{thm}
\noindent
{\em Proof.} As the proof is very similar to those of Theorem \ref{stab1}, we limit ourselves to 
pointing out the only differences. For the wellposedness, we use the estimate valid
for $s\ge 0$
\[
|| ( 1- \partial _x ^2 )^{-1} (au_x)_x||_s \le C ||u||_s.
\]
For $s=1$, we claim that \eqref{A52} is justified by noticing that for $u\in H^1(\T )$ 
\[
\langle -(au_x)_x ,  u \rangle _{H^{-1},H^1} = \langle au_x, u_x \rangle 
_{L^2,L^2}.
\]
Then the wellposedness statement is established as in Theorem \ref{stab1}.
Let us proceed to the asymptotic behavior. Pick any 
$u_0\in H^1(\T )$ and any $v_0\in H^1(\T )$ which is the weak limit in $H^1(\T )$ of a sequence $u(t_n)$ with 
$t_n\to \infty$ and $t_{n+1}-t_n\ge T$. Let us still denote by
$v$ the solution of \eqref{A2}. 
Equation \eqref{A37} has to be replaced by 
\begin{equation}
\label{A56}
\int_0^T\!\!\! \int_{\T } a(x) |v_x(x,t)|^2dxdt=0.
\end{equation}
To prove \eqref{A56}, we notice first that $u(t_n+\cdot )\to v$ in 
$C([0,T],H^s(\T ))$ for $s<1$, and that $||u(t_n+\cdot )||_{L^2(0,T;H^1(\T ))}
\le const$, so that, extracting a subsequence if needed, we have that
\[
u(t_n+\cdot ) \to v \qquad \text{ weakly in } L^2(0,T; H^1(\T )). 
\]
This yields, with \eqref{A52}, 
\[
\int_0^T\!\!\!\int_{\T} a(x)|v_x|^2dxdt 
\le \liminf_{n\to \infty}\int_{t_n}^{t_n+T}\!\!\!\int_{\T} a(x) |u_{nx}|^2dxdt=0.
\]
Therefore $v$ solves 
\begin{eqnarray}
&& v_t-v_{xxt} +v_x+vv_x =0,\qquad x\in\T,\ t\in (0,T),\label{A60}\\ 
&& v_x=0,\qquad x\in\omega, \ t\in (0,T), \label{A61}\\
&& v\in C([0,T]; H^s(\T ))\qquad \text{ for } s<1. \label{A62}
\end{eqnarray}
\eqref{A60} and \eqref{A61} yield $v_t=v_x=0$ in $\omega \times (0,T)$, hence
\[
v(x,t)=C\qquad \text{ for } (x,t)\in \omega\times (0,T)
\]
for some constant $C\in \R$. If the UCP Conjecture is still valid when $v$ is constant on the band $\omega \times (0,T)$, 
then 
$v\equiv C$. As $d[u(t)]/dt=0$, it follows that $C=[u_0]$ and 
that as $t\to \infty$, $u(t)\to [u_0]$ weakly in $H^1(\T)$ and strongly in $H^s(\T)$ for any $s<1$. \cqfd
\begin{rmk}
Similar results, but with convergences towards 0, hold for the system
\begin{eqnarray*}
&&u_t-u_{xxt}+u_x+uu_x-(a(x)u_x)_x=0,\\
&&u(0,t)=u(2\pi, t)=0,\\
&&u(x,0)=u_0(x)
\end{eqnarray*}
provided that $a(x)>0$ for $x\in (0,\varepsilon )\cup (2\pi -\varepsilon ,0)$ for some $\varepsilon >0$ and the UCP Conjecture holds. 
\end{rmk}
\subsection{Boundary stabilization of BBM}
In this section, we are concerned with the boundary stabilization of the BBM 
equation. We consider the following system
\begin{eqnarray}
&& u_t - u_{xxt} +u_x + uu_x =0, \qquad x\in (0,L),\ t\ge 0,
\label{B1}\\
&& u_{tx}(0,t)=\alpha u(0,t) + \frac{1}{3} u^2(0,t), \label{B2}\\
&& u_{tx}(L,t)=\beta u(L,t) + \frac{1}{3} u^2(L,t), \label{B3}\\
&& u(x,0)=u_0(x), \label{B4}
\end{eqnarray}
where $\alpha$ and $\beta$ are some real constants chosen so that 
\[
\alpha >\frac{1}{2} \quad \text{ and }\quad \beta < \frac{1}{2}\cdot \label{B5}
\]
Scaling in \eqref{B1} by $u$ yields (at least formally)
\begin{eqnarray}
\frac{d}{dt} \frac{1}{2} \int_0^L (|u|^2 + |u_x|^2)dx 
&=& \left[ uu_{tx}\right] _0^L -\left[ 
\frac{u^2}{2} +\frac{u^3}{3}\right] _0^L \nonumber\\
&=& ( \beta -\frac{1}{2}) |u(L,t)|^2 +(\frac{1}{2}-\alpha )|u(0,t)|^2,
\label{B6}
\end{eqnarray}
hence $||u(t)||_{H^1}$ is nonincreasing if $\alpha$ and $\beta$ 
fulfill \eqref{B5}. The wellposedness of \eqref{B1}-\eqref{B4}
and the asymptotic behavior are described in the following result.
\begin{thm}
Let $s\in (1/2,5/2)$ and $u_0\in H^s(0,L)$. Then there exist a time
$T>0$ and a unique solution $u\in C([0,T],H^s(0,L))$ of
\eqref{B1}-\eqref{B4}. Furthermore, if $s=1$, then $T$ may be taken arbitrarily
large, and if the UCP Conjecture holds, as $t\to \infty$
\begin{eqnarray}
&&u(t)\to 0 \quad \text{\rm weakly in } H^1(0,L) \label{B7}\\ 
&&u(t)\to 0 \quad \text{\rm strongly in } H^s(0,L) \text{ for all } s<1.\label{B8} 
\end{eqnarray}
\end{thm}
{\em Proof.} Let us begin with the well-posedness part. 
Pick any $u_0\in H^s(0,L)$ with $s>\frac{1}{2}$. Let $v=u_t$. 
Then $v$ solves the elliptic problem
\begin{eqnarray}
&&(1-\partial _x ^2)v=f, \quad x\in (0,L),   \label{B9}\\
&&v_x(0)=a,\ v_x(L)=b     \label{B10}
\end{eqnarray}
with $f:=-u_x-uu_x$, 
$a:=\alpha u(0,t) + \frac{1}{3} u^2(0,t)$,
$b:=\beta  u(L,t) + \frac{1}{3} u^2(L,t)$. Note that the solution 
$v$ of \eqref{B9}-\eqref{B10} may be written as
\[
v=w+g
\]
where 
$g(x)=ax+ \frac{b-a}{2L}x^2$ and 
$w=(1-\partial _x^2)_N^{-1}(f-(1-\partial ^2 _x)g)$ solves 
\begin{eqnarray*}
&&(1-\partial _x ^2)w=f-(1-\partial _x^2)g\\
&&w_x(0)=w_x(L)=0.
\end{eqnarray*}
Thus
\begin{equation}
u_t=v=-(1-\partial _x^2)_N^{-1} (u_x+uu_x) 
+ (1-(1-\partial _x^2)_N^{-1} (1-\partial_x^2))g.
\label{B15}
\end{equation}
We seek $u$ as a fixed point of the integral equation
\begin{eqnarray} 
u(t)={\Gamma }(u)(t)  &:=& u_0 + \int_0^t 
\bigg\{ 
-(1-\partial _x^2)_N^{-1}(u_x+uu_x)(\tau) \nonumber \\
&&\  \left. +(1-(1-\partial _x ^2)^{-1}_N (1-\partial _x^2))
\bigg[  [ \alpha u(0,\tau) +\frac{1}{3}u^2(0,\tau) ] x \right. \nonumber\\
&&   +(2L)^{-1} [  \beta u(L,\tau ) + \frac{1}{3} u^2(l,\tau ) 
-\alpha u(0,\tau ) -\frac{1}{3} u^2(0, \tau ) ] x^2 \bigg]\bigg\} d\tau . \quad\  \label{B16}
\end{eqnarray}

Note that ${\mathcal D}((1-\partial _x^2)_N^{\frac{s}{2}}) =H^s(0,L)$  for
$-1/2<s<3/2$. 
Let $R>0$ be given and let $B_R$ denote the closed ball in $C([0,T],H^s(0,L))$
of center $0$ and radius $R$. For $1/2<s<5/2$ and $u\in B_R,v\in B_R$, we have that
\begin{eqnarray*}
||{\Gamma }(u)(t)-{\Gamma }(v)(t)||_{H^s(0,L)} 
&\le& C\int_0^t (1+R) ||u(\tau ) -v(\tau )||_{H^s(0,L)} d\tau \\
&\le& CT(1+R) ||u-v||_{C([0,T],H^s(0,L))}  
\end{eqnarray*}
and
\[
||{\Gamma }(u)(t)||_{H^s(0,L)} \le ||u_0||_{H^s(0,L)} 
+CT (1+R) ||u||_{C([0,T],H^s(0,L))}\cdot
\]
Pick $R=2||u_0||_{H^s(0,L)}$ and $T=(2C(1+R))^{-1}$. Then
$\Gamma $ is a contraction in $B_R$, hence 
it admits a unique fixed point which solves \eqref{B16}, or 
\eqref{B15} and \eqref{B4}. 

Assume now that $s=1$. It follows from \eqref{B15} that 
$u_t\in C([0,T],H^2(0,L))$, hence we can scale by $u$ in \eqref{B1} to derive
\eqref{B6}. Thus $||u(t) ||_{H^1(0,L)}$ is nonincreasing, and $T$ may be taken as
large as desired. Let us turn our attention to the asymptotic behavior.
Let $(t_n)_{n\ge 0}$ be any sequence with $t_n\to \infty$ as $n\to \infty$.
Extracting a subsequence if needed, we can assume that $t_{n+1}-t_n\ge T$
for each $n$ and that, for some $v_0\in H^1(0,L)$,  $u(t_n)\to v_0$ weakly in $H^1(0,L)$ (hence strongly in $H^s(0,L)$ for $s<1$) as $n\to \infty$. The 
continuity of the flow map (which follows at once from the fact that the map 
$\Gamma$ is a contraction)  yields
\[
u(t_n+\cdot ) \to v(\cdot ) \qquad \text{\rm in }\ 
C([0,T],H^s(0,L))  \quad \text{\rm for } s<1,
\]
where $v$ denotes the solution of \eqref{B1}-\eqref{B4} 
issued from $v_0$ at $t=0$. Integrating \eqref{B6} on $[t_n,t_{n+1}]$ 
yields
\[
\frac{1}{2} ||u(t_{n+1} ) ||^2_{H^1(0,L)} -\frac{1}{2}||u(t_n) ||^2_{H^1(0,L)}
+(\frac{1}{2} -\beta ) \int_{t_n}^{t_{n+1}} |u(L,t)|^2dt
+(\alpha -\frac{1}{2})\int_{t_n}^{t_{n+1}} |u(0,t)|^2 dt. 
\]
Letting $n\to \infty$ yields
\begin{equation}
\label{B30}
(\frac{1}{2}-\beta ) \int_0^T |v(L,t)|^2 dt 
+(\alpha - \frac{1}{2} ) \int_0^T |v(0,t)|^2 dt =0. 
\end{equation}
Extending $v(x,t)$ by $0$ for $x\in \R \setminus (0,L)$  and $t\in (0,T)$,
we infer from \eqref{B1}, \eqref{B2}, \eqref{B3} (for $v$)  and 
\eqref{B30} that 
\[
v_t - v_{txx} + v_x +vv_x =0 \qquad \text{ for } x\in \R ,\ t\in (0,T)
\]
with
\[
v(x,t)=0 \qquad \text{ for } x\in \R \setminus (0,L),\ t\in (0,T).
\]
Since $v\in C([0,T],H^s(\R ))$ with $1/2 <s<1$, we infer from 
the UCP Conjecture (if true) that $v\equiv 0$, hence $v_0=0$. \cqfd

\end{document}